\documentclass{article}
\usepackage[utf8]{inputenc}
\usepackage[english]{babel}
\usepackage{latexsym}
\usepackage{amsfonts}
\usepackage{soul}
\usepackage{mathrsfs}
\usepackage{amsthm}
\usepackage[all]{xy}
\usepackage{mathtools}
\usepackage{color}
\usepackage{amssymb, mathrsfs, amsfonts, amsmath}
\usepackage{mathtools}
\usepackage{graphicx}
\usepackage[font=small,labelfont=bf]{caption}
\usepackage{epstopdf}
\usepackage[pdfpagelabels,hyperindex]{hyperref}
\usepackage{xcolor}
\usepackage{float}
\usepackage{pgfplots}
\usepackage{listings}
\usepackage{longtable}
\usepackage{mathrsfs}
\usepackage{amsbsy}
\usepackage{tikz}
\usepackage{array,tabularx,float, multicol,multirow,graphicx}
\usetikzlibrary{calc,intersections,arrows.meta,through,patterns}
\usetikzlibrary{decorations.pathmorphing}
\usepackage{pgfplots}
\pgfplotsset{compat=1.18}
\usetikzlibrary{calc,through,intersections,patterns,patterns.meta}
\usepackage{wasysym}
\usepackage{stackrel}
\pgfplotsset{compat=newest}
\usetikzlibrary{calc,through,intersections,patterns,patterns.meta}

\newtheorem{theorem}{Theorem}[section]
\newtheorem{thm}{Theorem}[section]
  
 \newtheorem{lem}[thm]{Lemma}
 \newtheorem{remark}[thm]{Remark}
 \newtheorem{prop}[thm]{Proposition}

 \newtheorem{definition}[thm]{Definition}
 \newcommand{\grad}{\nabla\,}
 \newcommand{\diver}{{{\rm{div}}}}
 \newcommand\raisepunct[1]{\,\mathpunct{\raisebox{0.5ex}{#1}}}
\usepackage{graphicx} 

\usepackage[english]{babel}

\usepackage[letterpaper,top=2cm,bottom=2cm,left=3cm,right=3cm,marginparwidth=1.75cm]{geometry}

\usepackage{amsmath}
\usepackage{graphicx}

\providecommand{\keywords}[1]{%
  \small\textbf{Keywords ---} #1
}

\title{$p$-Eigenvalue pinching sphere theorems}

\author{%
Paulo Henryque C. Silva \thanks{Universidade Federal do Ceará. Departamento de Matemática. Fortaleza - CE, Brazil. \texttt{paulohenryque@alu.ufc.br}}
}

\begin{document}
\maketitle
    
\begin{abstract}

\noindent In this paper, we establish two $p$-eigenvalue pinching sphere theorems,  for the \( p \)-Laplacian, $p>1$. The first result states that if the first non-zero $p$-eigenvalue of a closed Riemannian $n$-manifold with sectional curvature $K_{M}\geq 1$   is sufficiently close to the  first non-zero $p$-eigenvalue of $\mathbb{S}^{n}$  then $M$ is homeomorphic to $\mathbb{S}^{n}$. The second states that if the first non-zero $p$-eigenvalue of a closed Riemannian $n$-manifold with Ricci curvature ${\rm Ric}_{M}\geq (n-1)$ and injectivity radius ${\rm inj}_{M}\geq i_0>0$   is sufficiently close to the  first non-zero $p$-eigenvalue of $\mathbb{S}^{n}$  then $M$ is diffeomorphic to $\mathbb{S}^{n}$. Our results extend sphere theorems originally settled for the Laplacian by S. Croke~\cite{Croke1982} and G.P. Bessa~\cite{bessa} respectively. \\

\noindent\textbf{AMS Subject Classification:} 58C40, 58J50, 53C21, 35P30.

\noindent \keywords{Pinching Theorem; $p$‑eigenvalue; \( p \)-Laplacian; riemannian manifold.}

\end{abstract}

\section{Introduction} 
 There has been an increasing interest in the study of the operator $\triangle_{p}=\diver (\vert \grad u \vert^{p-2}\grad u )$,  $p>1$, called the $p$-Laplacian, on compact Riemannian manifolds with or without boundary and its relations with their geometries, particularly the relations regarding the spectrum of $\triangle_p$, which have attracted considerable attention in the recent literature (see \cite{BianchiPigolaSetti2020}, \cite{LimaMontenegroSantos2010}, \cite{ColboisProvenzano2021} e \cite{LiWang2021}). The spectrum of $\triangle_p$ is the set of $\lambda \in \mathbb{R}$ such that there exists $0\not \equiv u \in W^{1,p}_{0}(M)$ satisfying the following  equation \begin{equation} \label{eq1aa}\left\{\begin{array}{rl}
    \triangle_p u + \lambda \vert u\vert^{p-2}u=0  &  in \,\, \,\,M;\\
    u=0 & on \,\,\,\, \partial M,\,\,\, if\,\,\partial M \,\,\neq \emptyset.
\end{array}\right.
\end{equation} In weak sense, this is 
\begin{equation}\label{eqweak}
\int_{M} | \nabla u |^{p-2}\langle \grad u, \grad \eta\rangle \,d\mathscr{H}^{n}=\lambda \int_{M}\vert u\vert^{p-2} u \,\eta\, d\mathscr{H}^{n},   
\end{equation}
for all $\eta \in C^{\infty}_o(M)$. If $\partial M \neq \emptyset$ the 
problem \eqref{eq1aa} is called the Dirichlet eigenvalue problem and if $\partial M =\emptyset$, then the problem is called the closed eigenvalue problem.
Here, the Sobolev space \( W^{1,p}_0(M) \) is the completion of \( C^\infty_0(M) \) with respect to the norm
\[
\| u \|_{1,p} = \left( \int_M \left( |u|^p + | \nabla u |^p \right) \, d \mathscr{H}^{n} \right)^{1/p}.
\]

 In the Dirichlet eigenvalue problem, there exists a non-decreasing sequence of positive eigenvalues  $\lambda_{k,p}(M)\to \infty $ tending to infinity as $k \to \infty$, proved using the Ljusternik–Schnirelman principle see \cite{garzia-peral}, \cite{le} and references therein. The first Dirichlet eigenvalue (smallest) is a Ljusternik–Schnirelman eigenvalue, isolated and simple, meaning that there exists a $C^{1,\alpha}(M)$ function $u$ such that any other first  Dirichlet eigenfunction is a multiple of $u$, moreover, any first  Dirichlet eigenfunction does not change sign, see \cite{anane}, \cite{lindqvist}. The first Dirichlet eigenvalue of $M$ has a variational characterization given by 
 \begin{equation}\label{fristfrequence}
     \lambda_{1,p}^{\mathcal{D}}(M)=\inf\left\{\frac{ \displaystyle \int_{M} | \nabla u |^{p}d\mathscr{H}^{n}}{ \displaystyle \int_{M}|u|^{p}}d\mathscr{H}^{n}\colon 0 \not \equiv u \in W_{0}^{1,p}(M)\right\}. 
 \end{equation}

 For the closed eigenvalue problem, $\lambda=0$ is an eigenvalue corresponding to constant functions, hence, we consider as the first non-zero eigenvalue of \eqref{eq1aa}  the value given by 
 \begin{equation}\label{closeeigen}
   \lambda_{1,p}(M)=\left\{ \frac{\displaystyle \int_{M} | \nabla u |^{p}\,d\mathscr{H}^{n}}{ \displaystyle \int_{M}|u|^{p}} d\mathscr{H}^{n}\colon 0 \not \equiv u \in W^{1,p}(M);\,\, \int_{M} |u|^{p-2}u \, d\mathscr{H}^{n} =0 \right\}.   
 \end{equation}
  Since $\triangle_{2}=\triangle$  it is expected that many results that are valid for the Laplacian are also valid for the $p$-Laplacian. In \cite{matei}, Ana Matei and Hiroshi Takeuchi in  \cite{takeuchi} proved Cheng's eigenvalue comparison theorem and Faber-Krahn isoperimetric inequality for the $p$-Laplacian.  An L\^{e}  studied the spectrum of the $p$-Laplacian on bounded domains $\Omega$ with $\partial \Omega \neq \emptyset$ under various boundary conditions, see  \cite{le}.
  Regarding the closed eigenvalue problem, Matei \cite[Thm. 3.1]{matei} proved a Lichenorowicz-Obata version theorem for the $p$-Laplacian. More precisely, she proved the following result. \begin{theorem}[Matei-Lichenorowicz-Obata]\label{matei} Let $M$ be a closed Riemannian $n$-manifold with Ricci curvature ${\rm  Ric}_{M}\geq (n-1)$ and $1<p<\infty$ then \[\lambda_{1,p}(M) \geq \lambda_{1,p}(\mathbb{S}^{n}).\] With equality if and only if $M$ is isometric to $\mathbb{S}^{n}$.
    \end{theorem}In view of Theorem \eqref{matei} one  expects that a closed Riemannian manifold with  ${\rm  Ric}_{M}\geq (n-1)$, $1<p< \infty$, and first non-zero eigenvalue $C(n,p) \cdot \lambda_{1,p}(\mathbb{S}^{n}) \geq \lambda_{1,p}(M) $, for some suitable constant $C(n,p)>1$, is homeomorphic/diffeomorphic to the sphere. Our results essentially confirm this. 

 \begin{theorem}Let $M$ be a closed Riemannian $n$-manifold with sectional curvature $K_M\geq 1$. There exists a constant $C(n,  p)>1$ such that if $C(n,p)\cdot\lambda_{1,p}(\mathbb{S}^{n}) \geq \lambda_{1,p}(M)$ then $M$ is homeomorphic to $\mathbb{S}^{n}$.\label{thm2}
 \end{theorem}

\begin{theorem}\label{thm1}Let $M$ be a closed Riemannian $n$-manifold with ${\rm Ric}_{M}\geq(n-1)$ and injectivity radius ${\rm inj}_{M}\geq i_0>0$. There exists a constant $C(n, i_0, p)>1$ such that if $C(n,i_0,p)\cdot\lambda_{1,p}(\mathbb{S}^{n}) \geq \lambda_{1,p}(M)$ then $M$ is diffeomorphic to $\mathbb{S}^{n}$.
\end{theorem} 
 \begin{remark}Theorems \ref{thm2} and \ref{thm1} are natural versions of the $p$-Laplacian of the sphere theorems for Ricci curvature proved in  C. Croke~\cite{Croke1982} and~\cite{bessa}.  \end{remark}
 \begin{remark} The hypothesis ${\rm inj}_{M}\geq i_0>0$ cannot be weakened to $\mathscr{H}^{n}(M)\geq v>0$. There are
closed Riemannian manifolds $n$-manifolds $(M, g_{\epsilon})$, $n>4$  satisfying ${\rm Ric}_{(M,g_{\epsilon}}) \geq  (n-1)$, ${ \mathscr{H}_{ g_{\epsilon}}^{n}}(M)\geq v>0$ and ${\rm diam}(M, g_{\epsilon})\geq \pi - \epsilon$, for any $\epsilon >0$ which are not homotopy equivalent to $\mathbb{S}^{n}(1)$, see  \cite{Anderson}, \cite{otsu}.
  \end{remark}
  
\section{Preliminaries}

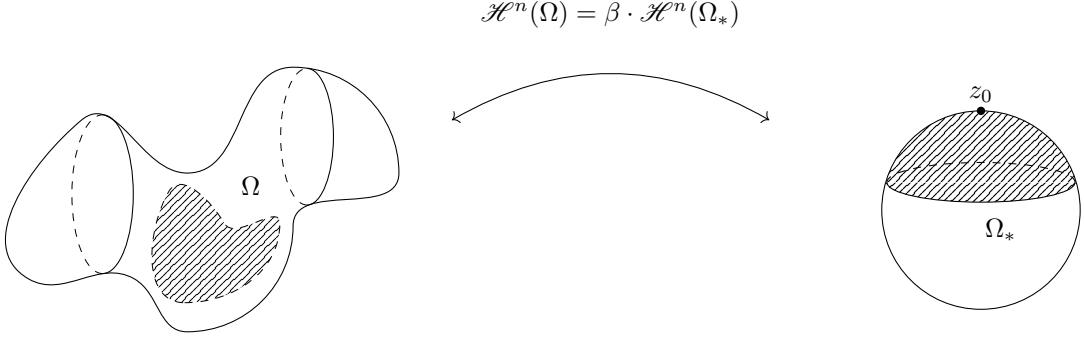
\begin{figure}[htb!]
    \centering

	\def\r{5}
	\begin{tikzpicture}[scale=0.70]
		
		\path[draw] (-1,-1) to[out=0,in=-90] (1,1) to[out=90,in=-90] (3,2) to[out=90,in=0] (1,4) to[out=180,in=0] (-1,2) to[out=180,in=40] (-3,3) to[out=220,in=150] (-4,0) to[out=-30,in=150] (-2,0) to[out=-30,in=180] (-1,-1);
		
		\path[draw,name path=arc1] (-2.6,0.1) arc (-90:90:0.58cm and 1.5 cm);
		\path[draw,name path=arc1,dashed,very thin] (-2.6,0.1) arc (270:90:0.58cm and 1.5 cm);
		
		\begin{scope}[xshift=3.5cm,yshift=1.31cm,scale=0.86]
			\path[draw,name path=arc1] (-2.6,0.1) arc (-90:90:0.58cm and 1.5 cm);
			\path[draw,name path=arc1,dashed,very thin] (-2.6,0.1) arc (270:90:0.58cm and 1.5 cm);
		\end{scope}
		
		\begin{scope}[scale=1.5,xshift=-0.2cm,yshift=-0.3cm]
			\path[draw,dashed,very thin,pattern={Lines[angle=45,distance={2pt}]}] (-0.5,0) to[out=0,in=-90] (0.7,1) to[out=90,in=-60] (0,1) to[out=120,in=60] (-0.7,1.4) to[out=240,in=180] (-0.5,0);
		\end{scope}
		\draw (0.2,1.4) node[above] {$\Omega$};

		\begin{scope}[xshift=14cm,yshift=1.3cm,scale=0.75]
			\draw (0,0) circle (2.5cm);
			\path[draw,name path=arc1,dashed] (2.37,0.7) arc (0:180:2.37cm and 0.5 cm);
			\path[draw,name path=arc1] (2.37,0.7) arc (360:180:2.37cm and 0.5 cm);
			
			\begin{scope}
				\clip (0,0) circle (2.5 cm);
				\path[pattern={Lines[angle=45,distance={2pt}]}] (2.4,0.7) arc (360:180:2.4cm and 0.5 cm)--(2.4,0.7) arc (0:180:2.4cm and 1.79 cm);
			\end{scope}

			\node (t) at (0.5,-0.5) {$\Omega_{\ast}$};
			\node[above] (t) at (0,2.5) {$z_{0}$};
			\fill (0,2.5) circle (3pt);
			
		\end{scope}
		\draw[<->] (4,3) to[out=30,in=150] (10,3);
		\node[] (tt) at (7,5)  {\normalsize{$\mathscr{H}^{n}(\Omega) = \beta \cdot \mathscr{H}^{n}(\Omega_*)$}};
			
	\end{tikzpicture}
    
    \caption{Schwarz symmetrization of $\Omega$ }
    \label{figSC}
\end{figure}

Let \( M \) be an \( n \)-dimensional Riemannian manifold such that \( \mathrm{Ric}_M \geq (n - 1) \). Initially, we denote by \( \beta \) the following number.
  \[ \displaystyle \beta \coloneqq \frac{\mathscr{H}^{n}(M)}{\mathscr{H}^{n}(\mathbb{S}^n)}\raisepunct{,}
\]
where \( \mathscr{H}^{n} \) denotes the \( n \)-dimensional Hausdorff measure. Let \( \Omega \subset M \) be an open subset with a smooth boundary. The \emph{Schwarz symmetrization} of \( \Omega \) consists in the following construction: define \( \Omega_* \subset \mathbb{S}^n \) as the geodesic ball centered at the north pole \(z_0\) of the unit sphere \( \mathbb{S}^n \) such that its total volume is given by the following expression:
\begin{equation}\label{symm}
    \mathscr{H}^{n}(\Omega) = \beta \cdot \mathscr{H}^{n}(\Omega_*).
\end{equation}
Geometrically, Schwarz symmetrization can be described as illustrated in the figure \ref{figSC}.

Under these conditions \textit{Gromov's isoperimetric inequality} \cite{gromov}, gives us
\begin{equation}\label{gromov}
\mathscr{H}^{n-1}(\partial \Omega) \geq  \beta\cdot \mathscr{H}^{n-1}(\partial \Omega_*).    
\end{equation}
In \cite{Croke1982}, Croke establishes a refined version of Gromov's isoperimetric inequality under the additional assumption that \( \operatorname{diam}(M)\leq D < \pi \). More precisely, he proves the following result:
\begin{lem}[Croke]\label{lemacroke}
Let \( M \) be a Riemannian $n$-manifold with Ricci curvature \( \mathrm{Ric}_M \geq n - 1 \), and  diameter  \( \operatorname{diam}(M) \leq D < \pi \). Then, there exists a constant \( \overline{C}(n, D) > 1 \) such that, for any domain \( \Omega \subset M \) with smooth boundary \( \partial \Omega \), we have
\[
\mathscr{H}^{n-1}(\partial \Omega) \geq \overline{C}(n, D) \cdot \beta \cdot \mathscr{H}^{n-1}(\partial \Omega_*),
\]
where \( \beta \) and \( \Omega_* \) are as previously defined.
\end{lem}
This lemma will be a key step in proving Theorems~\ref{thm1} and~\ref{thm2}. In addition to Lemma~\ref{lemacrokee}, two other tools play a central role in our approach. Although both are well known in the classical literature, we will provide a brief introduction to each tool for sake of completeness. The first is the coarea formula, which will be discussed below; the second is symmetrization, which will also be briefly addressed later.
Let \( M \) be a smooth Riemannian manifold  and \( f\colon M \rightarrow \mathbb{R} \)  be a \( C^1 \) function without critical points.  Then, for any measurable function \( \psi: M \rightarrow \mathbb{R} \), we have
\begin{equation}\label{coareaform}
\int_M \psi \, d\mathscr{H}^{n} = \int_{\mathbb{R}} \left( \int_{\{f = t\}} \frac{\psi}{| \nabla f |} \, d\mathscr{H}^{n-1} \right) dt.  
\end{equation}
In particular, by taking \( \psi = 1 \), we obtain
\begin{equation}\label{coareaformvol}
\mathscr{H}^{n}(M) = \int_{\mathbb{R}} \left( \int_{\{f = t\}} \frac{1}{| \nabla f |} \, d\mathscr{H}^{n-1} \right) dt.    
\end{equation} For more details see \cite[Section 6.2]{SimonGMT}.

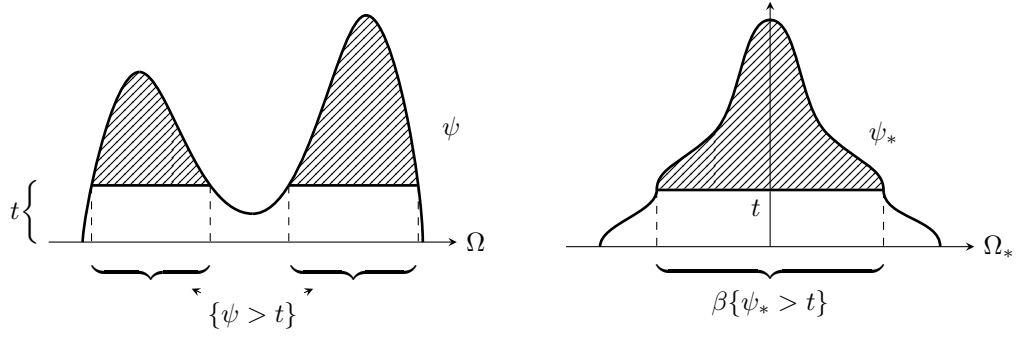
\begin{figure}[H]
    
    \caption{Deacreasing symmetric rearrangement}
   \begin{minipage}[h]{6cm}
   \hspace{0.20\textwidth}
	\begin{tikzpicture}[scale=0.75, >=stealth]
		\def\v{1.0}
		\def\p{3.6}
		
		\coordinate (o) at (0,0);
		\draw[->] ($(o)+(-1*\p,0)$) -- (\p,0) node[right] {$ \Omega $};
		\node (t1) at (3.5,2) {$ \psi $};
		
		\coordinate (p1) at (0,0.5);
		\coordinate (p2) at (1,2);
		\coordinate (p3) at (2,4);
		\coordinate (p4) at (3,0);
		
		\coordinate (p5) at (-1,2);
		\coordinate (p6) at (-2,3);
		\coordinate (p7) at (-3,0);
		
		\coordinate (c1) at ($(p1)+(1*\v,0)$);
		
		\coordinate (c2) at ($(p2)+0.5*(-1,1)$);
		\coordinate (c3) at ($(p2)+0.5*(1,-1)$);
		
		\coordinate (c4) at ($(p3)+(-0.5*\v,0)$);
		\coordinate (c5) at ($(p3)+(0.5*\v,0)$);
		
		\coordinate (c6) at ($(p4)+(0,\v)$);
		\coordinate (c7) at ($(p4)+(\v,0)$);
		
		\coordinate (c8) at ($(p1)+(-1*\v,0)$);
		
		\coordinate (c9) at ($(p5)+(\v,0)$);
		\coordinate (c10) at ($(p5)+0.5*(-1,-1)$);
		
		\coordinate (c11) at ($(p6)+(0.5*\v,0)$);
		\coordinate (c12) at ($(p6)+(-0.5*\v,0)$);
		
		\coordinate (c13) at ($(p7)+(0,0.4*\v)$);

		\begin{scope}
			\path[clip,name path global=curva] (p4) .. controls (c6) and (c5) .. (p3) .. controls (c4) and (c1) .. (p1) .. controls (c8) and (c11) .. (p6) .. controls (c12) and (c13) .. (p7) -- (p4);
			
			\path[fill,black!10,pattern=north east lines, draw,line width=1pt,name path global=rec1] (-3,1) rectangle (3,4);
			
			\draw[line width=1pt] (-3,1) rectangle (3,4);

		\end{scope}
		
		\path[name intersections={of=curva and rec1}];
		\coordinate (i2) at (intersection-2);
		\coordinate (i3) at (intersection-3);
		\coordinate (i4) at (intersection-4);
		\coordinate (i5) at (intersection-5);

		\draw[dashed] (i2) --  (i2|-o);
		\draw[dashed] (i3) --  (i3|-o);
		\draw[dashed] (i4) --  (i4|-o);
		\draw[dashed] (i5) --  (i5|-o);
		
		\node[below=2pt]  (t3) at ($(i2)!0.5!(i3)+(0,-1.1)$) {$ \underbrace{\hspace{1.65cm}}$};
		
		\node[below=2pt]  (t4) at ($(p1)+(0,-1.3)$) {$ \{\psi>t\}$};
		
		\node[below=2pt]  (t5) at ($(i5)!0.5!(i4)+(0,-1.1)$) {$ \underbrace{\hspace{1.55cm}}$};
		
		\draw[->] (t4) -- (t3);
		
		\draw[->] (t4) -- (t5);

		\node (t6) at ($(i5)!0.47!(p7)+(-1.1,0)$) {$t\bigg\{$};
		
		\draw[line width=1pt] (p1) .. controls (c1) and (c4) .. (p3);
		
		\draw[line width=1pt] (p3) .. controls (c5) and (c6) .. (p4);
		
		\draw[line width=1pt] (p1) .. controls (c8) and (c11) .. (p6);
		
		\draw[line width=1pt] (p6) .. controls (c12) and (c13) .. (p7);
	\end{tikzpicture}
\end{minipage}\hfill	
\begin{minipage}[h]{7.4cm}
\vspace{-0.045\textwidth}
\centering
	\begin{tikzpicture}[scale=0.75, >=stealth]
		\def\v{1}
		\def\p{3.6}
		
		\coordinate (o) at (0,0);
		\draw[->] ($(o)+(-1*\p,0)$) -- (\p,0) node[right] {$ \Omega_* $};
		\draw[->] (o) --  ($(o)+(0,1.2*\p)$);
		\node (t1) at (2,2) {$ \psi_{*} $};
		\node[below left] (t2) at (0,1) {$ t $};
		\node[below=2pt]  (t3) at (0,0) {$ \underbrace{\hspace{3.0cm}}$};
		
		\node[below=2pt]  (t5) at (0,-0.5) {$\beta \{\psi_{*}>t\}$};
		
	
		\coordinate (p1) at (0,4);
		\coordinate (p2) at (1,2);
		\coordinate (p3) at (2,1);
		\coordinate (p4) at (3,0);
		
		\coordinate (p5) at (-1,2);
		\coordinate (p6) at (-2,1);
		\coordinate (p7) at (-3,0);

		\coordinate (c1) at ($(p1)+(0.5*\v,0)$);
		
		\coordinate (c2) at ($(p2)+0.5*(-1,1)$);
		\coordinate (c3) at ($(p2)+0.5*(1,-1)$);
		
		\coordinate (c4) at ($(p3)+(0,0.4*\v)$);
		\coordinate (c5) at ($(p3)+(0,-0.5*\v)$);
		
		\coordinate (c6) at ($(p4)+(0,0.5*\v)$);
		\coordinate (c7) at ($(p4)+(\v,0)$);
		
		\coordinate (c8) at ($(p1)+(-0.5*\v,0)$);
		
		\coordinate (c9) at ($(p5)+0.5*(1,1)$);
		\coordinate (c10) at ($(p5)+0.5*(-1,-1)$);
		
		\coordinate (c11) at ($(p6)+(0,0.4*\v)$);
		\coordinate (c12) at ($(p6)+(0,-0.4*\v)$);
		
		\coordinate (c13) at ($(p7)+(0,0.4*\v)$);

		\draw[dashed] (p3) --  (p3|-o);
		\draw[dashed] (p6) --  (p6|-o);
		
		

		\begin{scope}
			\clip (p1) .. controls (c1) and (c2) .. (p2) .. controls (c3) and (c4) .. (p3) .. controls (c5) and (c6) .. (p4) -- (p7)  .. controls (c13) and (c12) .. (p6)  .. controls (c11) and (c10) .. (p5) .. controls (c9) and (c8) .. (p1);
			
			\fill[black!10,pattern=north east lines, draw,line width=1pt] (-2,1) rectangle (2,4);
			
			\draw[line width=1pt] (-2,1) rectangle (2,4);
			
		\end{scope}
		
		\draw[line width=1pt] (p1) .. controls (c1) and (c2) .. (p2);
		
		\draw[line width=1pt] (p2) .. controls (c3) and (c4) .. (p3);
		
		\draw[line width=1pt] (p3) .. controls (c5) and (c6) .. (p4);
		
		\draw[line width=1pt] (p1) .. controls (c8) and (c9) .. (p5);
		
		\draw[line width=1pt] (p5) .. controls (c10) and (c11) .. (p6);
		
		\draw[line width=1pt] (p6) .. controls (c12) and (c13) .. (p7);
	\end{tikzpicture}
\end{minipage}	 
    \label{fig1}
\end{figure}

\begin{definition}\label{definsr}
Let \( \psi \colon M \to [0, \infty) \) be a measurable function and $\Omega(t)=\{x\in M\colon \psi(x)>t\}$. Let $\Omega_{*}(t) \subset \mathbb{S}^{n}$ be the geodesic ball centered at the north pole, so that
\begin{equation}\label{defsym}
    \mathscr{H}^{n}(\Omega(t))=\beta \cdot \mathscr{H}^{n}(\Omega_{*}(t)) .
\end{equation}
We define \( \psi_* \colon \mathbb{S}^n \to [0, \infty) \) to be the function satisfying the following two properties.
\begin{enumerate}
\item $\psi_{*}$ is radial and decreasing;
  \item $\psi_{*}=t$ on the boundary of $\partial \Omega_{*}({t})$. 
 \end{enumerate} 
\end{definition}
Based on the properties discussed above, we can now provide a geometric interpretation of the symmetric decreasing rearrangement, as illustrated in Figure \ref{fig1}, where \( \Omega \) and \( \Omega_* \) are defined as in equation \ref{symm}. 

Due to the nature of our problem, the functions \( u \) that we aim to symmetrize belong to the Sobolev space \( W^{1,p}(\Omega) \), where \( \Omega \) is a bounded domain. To ensure that the symmetrized function \( u_{*} \) also belongs to the corresponding space \( W^{1,p}(\Omega_{*}) \), we rely on results concerning the \( p \)-Dirichlet energy, which assert that the energy decreases under symmetrization. 

The version adopted here is a refinement of the classical result due to Hildén~\cite{Hilden}.
\begin{prop}[P\'{o}lya–Szeg\"{o}-type inequality]\label{prop:PolyaSzego}
Let \( f: \Omega \rightarrow \mathbb{R} \) be a non-negative function in \( W^{1,p}(\Omega) \) such that \( f = 0 \) on \( \partial \Omega \). Then, for any \( 1 \leq p < \infty \), it holds that
\begin{equation}\label{PolyaSzego}
\int_{\Omega} |\nabla f |^p \, d\mathscr{H}^{n}\geq \beta \cdot \int_{\Omega_*} |\nabla f_*|^p \, d\mathscr{H}^{n}. 
\end{equation}
Where \( f_* \) and \( \Omega_* \) denote the symmetrizations of \( f \) and \( \Omega \), respectively.
\end{prop}
A proof of this result can be found in \cite[Cor.~3.22]{Baernstein}. Inequality  \eqref{PolyaSzego} shows that the symmetrization of  \( f \in W^{1,p}(\Omega) \) belongs to the corresponding Sobolev space; this is \(f_* \in W^{1,p}(\Omega_*) \).
The following lemma is the extension of Croke's Lemma \cite{Croke1982} to the $p$-Laplacian $\triangle_p$.

\begin{lem}\label{lemacrokee}
Suppose that \( \mathrm{Ric}_M \geq n - 1 \), \( \operatorname{diam}(M) \leq D < \pi \), and \( 1 < p < \infty \). Then there exists a constant \( C(n, p, D) > 1 \) such that
\begin{equation}\label{ineprinc}
    \lambda_{1,p}(M) \geq C(n, p, D) \cdot \lambda_{1,p}(\mathbb{S}^n).
\end{equation}
\end{lem}

\noindent \begin{proof}Let $\Omega \subset M$ be an open subset with smooth boundary $\partial \Omega.$ As in (\ref{fristfrequence}), the first Dirichlet eigenvalue $\lambda_{1,p}^{\mathcal{D}}(\Omega)$  of $\Omega $ is given by
\[
\lambda_{1,p}^{\mathcal{D}}(\Omega) = \inf\left\{ \frac{ \displaystyle \int_{\Omega} |\nabla \varphi |^p \, d\mathscr{H}^{n}}{ \displaystyle \int_{\Omega} |\varphi|^p \, d\mathscr{H}^{n}} : \varphi \in W_0^{1,p}(\Omega) \setminus \{0\} \right\}.
\]
This infimum is attained by a function \( u \in W_0^{1,p}(\Omega) \); see~\cite[Lem.~5]{lindqvist92}. Let \( u \in W_0^{1,p}(\Omega) \) be such a minimizer. By standard regularity theory, it follows that \( u \in C^{1,\alpha}_{\mathrm{loc}}(\Omega) \); see~\cite[Thm.~4.5]{le}. Moreover, the first Dirichlet eigenvalue is simple in any bounded domain~\cite[Thm.~6]{lindqvist} and is uniquely characterized by the fact that its associated eigenfunction does not vanish in the interior of \( \Omega \); see~\cite[Lem.~5]{lindqvist}. Hence, we may assume without loss of generality that \( u > 0 \) in \( \Omega \). Furthermore, by classical density arguments, we may also assume that \( u \) is a Morse function; see~\cite[Cor.~A.4]{le}.

For each \( t > 0 \), set
\[
\Omega(t) := \{ x \in \Omega : u(x) > t \},
\quad
\Gamma(t) := \{ x \in \Omega : u(x) = t \},
\] and $\Omega_{*}(t)\subset \mathbb{S}^{n}$ be a geodesic ball with center at the north pole such that 
\begin{equation}\label{symmm}
\mathscr{H}^{n}(\Omega(t)) = \beta \cdot \mathscr{H}^{n}(\Omega_*(t)),    
\end{equation}
with \( \beta= \mathscr{H}^{n}(M)/\mathscr{H}^{n}(\mathbb{S}^n) \). Let $u_{*}\colon \Omega (0)\to \mathbb{R}$, defined by 
\begin{equation}\label{eq:levelsets}
u_{*}\big|_{\Gamma_{*}(t)} = t = u\big|_{\Gamma(t)} .
\end{equation}
Where \( \Gamma_*(t) = \partial (\Omega_*(t)) \) be the decreasing symmetric rearrangement of $u$.  By construction, we have \( u_*\big|_{\partial \Omega_*} = 0 \).  Applying the coarea formula \eqref{coareaformvol}, we obtain
\begin{equation}\label{coareaomega}
    \mathscr{H}^{n}(\Omega(t)) = \int_{t}^{\infty} \left( \int_{\Gamma(s)} \frac{1}{ | \nabla u |} \, d\mathscr{H}^{n-1} \right) ds.   
\end{equation}
Differentiating the equation above with respect to $t$, it follows that
\begin{eqnarray}\label{eq1.10} 
   \dfrac{d}{dt}\mathscr{H}^{n}(\Omega(t)) &=&   \dfrac{d}{dt} \int_{t}^{\infty} \left( \int_{\Gamma(s)}\dfrac{1}{| \nabla u |} d\mathscr{H}^{n-1} \right) ds \\ \nonumber
  &=& - \int_{\Gamma(t)} \dfrac{1}{| \nabla u |} \ d\mathscr{H}^{n-1}.
\end{eqnarray}
Applying the coarea formula \eqref{coareaform} once again, together with the previous result \eqref{eq1.10} and the identities \eqref{symmm}, \eqref{eq:levelsets}, and \eqref{coareaomega}, we find that
\begin{eqnarray}\label{eq1.14a} \nonumber
\int_{\Omega(t)} u^p \, d\mathscr{H}^{n} &=&  \int_t^{\infty} \left( \int_{\Gamma(s)}\dfrac{u^p}{| \nabla u |} d\mathscr{H}^{n-1} \right) ds \\  \nonumber
&=& \int_t^{\infty} s^p \left( \int_{\Gamma(s)}\dfrac{1}{| \nabla u |} d\mathscr{H}^{n-1} \right) ds \\ 
&=& -  \int_t^{\infty} s^p \left( \dfrac{d}{ds}\mathscr{H}^{n}(\Omega(s) ) \right)  ds \\  \nonumber
&=& -\beta \cdot \int_t^{\infty} s^p \left( \dfrac{d}{ds}\mathscr{H}^{n}(\Omega_*(s) ) \right)  ds \\
&=&  \beta  \cdot \int_t^{\infty}\left( \int_{\Gamma(s)}\dfrac{u_*^p}{|\nabla u_* |} d\mathscr{H}^{n-1} \right) ds \nonumber  \\
&=& \beta \cdot  \int_{\Omega_*(t)} u_*^p \, d\mathscr{H}^{n}. \nonumber 
\end{eqnarray}
By the coarea formula once again, it follows that
\begin{equation}\label{eq1.11a}
     \dfrac{d}{dt}\int_{\Omega(t)} | \nabla u |^p \, d\mathscr{H}^{n} = - \int_{\Gamma(t)} | \nabla u |^{p-1} \, d\mathscr{H}^{n-1}.
\end{equation}
 Given the symmetric decreasing rearrangement $u_*$, we can show that 
\begin{equation}\label{symmmfor}
     \dfrac{d}{dt}\mathscr{H}^{n}(\Omega_*(t)) = -  \int_{\Gamma_*(t)} \dfrac{1}{ | \nabla u_*|} d\mathscr{H}^{n-1},
\end{equation}
and also that,
\begin{equation}\label{eq1.12a}
     \dfrac{d}{dt}\int_{\Omega_*(t)} |\nabla u_* |^p \, d\mathscr{H}^{n} = - \int_{\Gamma_*(t)} |\nabla u_*|^{p-1} \, d\mathscr{H}^{n-1}.
\end{equation}
By virtue of the construction of \(\Omega_*(t)\) given in \eqref{symmm}
and the identities in \eqref{eq1.10} and \eqref{symmmfor}, we see that
\begin{equation}\label{eqsym}
\begin{split}
\int_{\Gamma(t)} \frac{1}{| \nabla u |} \, d\mathscr H^{n-1}
&= -\,\mathscr H^n\bigl(\Omega(t)\bigr) \\
&= -\,\beta\, \cdot \mathscr H^n\bigl(\Omega_*(t)\bigr) \\
&= \beta\, \cdot
\int_{\Gamma_*(t)} \frac{1}{| \nabla u_*|} \, d\mathscr H^{n-1}.
\end{split}
\end{equation}

Using Hölder’s inequality, we obtain:
\begin{align*}
\int_{\Gamma(t)} d\mathscr{H}^{n-1} 
&\leq \left( \int_{\Gamma(t)} | \nabla u |^{p-1} \, d\mathscr{H}^{n-1} \right)^{1/p}
     \left( \int_{\Gamma(t)} \dfrac{1}{ | \nabla u |}  \, d\mathscr{H}^{n-1} \right)^{(p-1)/p}.
\end{align*}
As a result, combining the above inequality with \eqref{eq1.11a}, we conclude that
\begin{eqnarray}\label{intine}
   -\dfrac{d}{dt}\int_{\Omega(t)} | \nabla u |^p \, d\mathscr{H}^{n} & \geq &\frac{ \displaystyle \left( \int_{\Gamma(t)} d\mathscr{H}^{n-1} \right)^p}{ \displaystyle\left( \int_{\Gamma(t)} \dfrac{1}{| \nabla u |}  \, d\mathscr{H}^{n-1} \right)^{p-1}} \ \cdot
   \end{eqnarray}

According to Croke's Lemma \ref{lemacroke}, there exists a constant \( \overline{C}(n, D) > 1 \) such that for any domain \( \Omega \subset M \) with  smooth boundary, the following Gromov-type inequality holds:
\begin{equation}\label{crokelema}
    \int_{\Gamma(t)} d\mathscr{H}^{n-1}  \geq \overline{C}(n,D) \cdot \beta \cdot \int_{\Gamma_{*}(t)} d\mathscr{H}^{n-1}.
\end{equation} Thus, \eqref{intine} combined with \eqref{crokelema} and \eqref{eqsym}   we have that
\begin{equation}\label{eq4.3a}
   -\dfrac{d}{dt}\int_{\Omega(t)} | \nabla u |^p \, d\mathscr{H}^{n} \geq \dfrac{ \displaystyle \overline{C}(n,D)^p \cdot\beta^p}{\beta^{p-1}}\cdot \frac{ \displaystyle   \left( \int_{\Gamma_*(t)} d\mathscr{H}^{n-1} \right)^p}{ \displaystyle \left(  \int_{\Gamma_*(t)}  \dfrac{1}{| \nabla u_*|}  d\mathscr{H}^{n-1} \right)^{p-1}} \ \cdot
\end{equation}
On the other hand, due to the fact that $u_{*}$ is radially symmetric, it follows that
\begin{equation}\label{eq4.4a}
\int_{\Gamma_*(t)} d\mathscr{H}^{n} 
= \left( \int_{\Gamma_*(t)} | \nabla u_*|^{p-1} \, d\mathscr{H}^{n} \right)^{1/p}
     \left( \int_{\Gamma_*(t)} \dfrac{1}{| \nabla u_*|}  \, d\mathscr{H}^{n} \right)^{(p-1)/p}.
\end{equation}
Indeed, it follows by H\"older's inequality that  
\begin{eqnarray*}
    \int_{\Gamma_*(t)} d\mathscr{H}^{n-1} &\leq& \left( \int_{\Gamma_*(t)} | \nabla u_*|^{p-1} \, d\mathscr{H}^{n-1} \right)^{1/p}
     \left( \int_{\Gamma_*(t)} \dfrac{1}{| \nabla u_*|}  \, d\mathscr{H}^{n-1} \right)^{(p-1)/p} \\
     &=& \left( \left| \dfrac{du_*}{dt} \right|^{p-1} | \nabla t |^{p-1} \right)^{1/p} \left( \left| \dfrac{du_*}{dt} \right|^{-1}\ | \nabla t |^{-1} \right)^{(p-1)/p} \left( \int_{\Gamma_*(t)} d\mathscr{H}^{n-1} \right) \\
     &=&  \int_{\Gamma_*(t)} d\mathscr{H}^{n-1}.
\end{eqnarray*} This shows \eqref{eq4.4a}.
Therefore, inequality \eqref{eq4.3a} becomes \begin{eqnarray*}
    -\dfrac{d}{dt}\int_{\Omega(t)} | \nabla u |^p \, d\mathscr{H}^{n} \geq \overline{C}(n,D)^p\cdot  \beta \cdot \left( \int_{\Omega_*(t)} | \nabla u_*|^{p-1} \, d\mathscr{H}^{n-1} \right),
\end{eqnarray*} and by \eqref{eq1.12a} we have
\begin{eqnarray*}
    -\dfrac{d}{dt}\int_{\Omega(t)} | \nabla u |^p \, d\mathscr{H}^{n} &\geq& \overline{C}(n,D)^p\cdot  \beta \cdot \left( -  \dfrac{d}{dt}\int_{\Omega_{*}(t)} | \nabla u_*|^p \, d\mathscr{H}^{n}\right).
\end{eqnarray*}  
Integrating both sides from $t$ to $\infty$, we obtain:
\begin{equation}\label{gradine}
    \int_{\Omega(t)} | \nabla u |^p d\mathscr{H}^{n} \geq \overline{C}(n,D)^p \cdot\beta \cdot \int_{\Omega_*(t)}  | \nabla u_*|^p  d\mathscr{H}^{n}.
\end{equation} 
By combining the equality \eqref{eq1.14a} with the preceding inequality \eqref{gradine}, we conclude  
\begin{eqnarray*}\lambda_{1,p}^{\mathcal{D}}(\Omega)&=&\frac{\displaystyle \int_{\Omega(0)} | \nabla u |^p \,d\mathscr{H}^{n}}{\displaystyle \int_{\Omega(0)}u^{p}\,d\mathscr{H}^{n} } \nonumber\\
&\geq &  \overline{C}(n,D)^p \cdot\beta \cdot \frac{\displaystyle \int_{\Omega_*(0)}  | \nabla u_*|^p  d\mathscr{H}^{n}}{\displaystyle \int_{\Omega(0)}u^{p}\,d\mathscr{H}^{n} } \nonumber \\
&=&\overline{C}(n,D)^p \cdot\beta \cdot \frac{\displaystyle \int_{\Omega_*(0)}  | \nabla u_*|^p  d\mathscr{H}^{n}}{\displaystyle \beta\cdot\int_{\Omega_{*}(0)}u_{*}^{p}\,d\mathscr{H}^{n} } \nonumber\\
&\geq & \overline{C}(n,D)^p \cdot\lambda_{1,p}^{\mathcal{D}}(\Omega_{*}).
    \end{eqnarray*}

    The first non-zero eigenfunction of the $p$-Laplacian on $M$ has exactly two nodal domains $\Omega$ and $\Omega^{c}$ where $$\lambda_{1,p}(M)=\lambda_{1,p}^{\mathcal{D}}(\Omega)=\lambda_{1,p}^{\mathcal{D}}(\Omega^{c}),$$   see  \cite[Lemma 3.2]{matei}, \cite{Veron}. Suppose, without loss of generality, that $$\mathscr{H}^{n}(\Omega)\leq  \mathscr{H}^{n}(M)/2 \leq \mathscr{H}^{n}(\mathbb{S}^{n})/2.$$ Since $\mathscr{H}^{n}(\Omega)=\beta \cdot \mathscr{H}^{n}(\Omega_{*})$, we have that $\mathscr{H}^{n}(\Omega_{*}) \leq \mathscr{H}^{n}(\mathbb{S}^{n})/2$. This means that $\Omega_{*}\subset \mathbb{S}^{n}_{+}$, the north hemisphere of $\mathbb{S}^{n}$. Therefore \[\lambda_{1,p}(M)=\lambda_{1,p}^{\mathcal{D}}(\Omega)\geq \overline{C}(n,D)^{\,p}\cdot \lambda^{\mathcal{D}}_{1,p}(\Omega_{*})\geq \overline{C}(n,D)^{\,p}\cdot\lambda^{\mathcal{D}}_{1,p}(\mathbb{S}^{n}_{+})=\overline{C}(n,D)^{\,p}\cdot\lambda_{1,p}(\mathbb{S}^{n}). \] Since $\lambda_{1,p}(\mathbb{S}^{n}_{+})=\lambda_{1,p}(\mathbb{S}^{n})$, see \cite[Cor. 3.1]{matei}.
Setting $C(n,p,D)=\overline{C}(n,D)^{\,p}$ we have the Lemma \ref{lemacrokee}.
\end{proof}

\section{Proof of Results}To prove Theorems \ref{thm1} and \ref{thm2} we need the following diameter sphere theorems. Theorem~\ref{Grove-Shio} can be found in \cite{GS}, while Theorem~\ref{G.P. Bessa} is given in \cite[Thm.~1.1]{bessa}.

\begin{theorem}[Grove–Shiohama Sphere Theorem \cite{GS}]\label{Grove-Shio}
Let \( M \) be a closed Riemannian \( n \)-manifold with sectional curvature \( K_M \geq 1 \). If \( \operatorname{diam}(M) \geq \pi / 2 \), then \( M \) is homeomorphic to \( \mathbb{S}^n \).
\end{theorem}

\begin{theorem}[Bessa] \label{G.P. Bessa}
Given an integer \( n \geq 2 \) and a constant \( i_0 > 0 \), there exists \( \epsilon = \epsilon(n, i_0) > 0 \) such that if \( M \) is an \( n \)-dimensional Riemannian manifold satisfying
\[
\operatorname{Ric}_{M} \geq n - 1, \quad \operatorname{inj}_{M} \geq i_0, \quad \operatorname{diam}(M) \geq \pi - \epsilon,
\]
then \( M \) is diffeomorphic to the standard sphere \( \mathbb{S}^n \) endowed with the canonical metric.
\end{theorem}
\subsection{Proof of Theorems \ref{thm2} and \ref{thm1}}

    If $M$ is a closed Riemannian $n$-manifold with Ricci curvature ${\rm Ric}_{M}\geq (n-1)$ then ${\rm diam}(M)\leq \pi$. By Lemma \ref{lemacrokee} if $\lambda_{1,p}(M) \leq C(n,p,D) \lambda_{1,p}(\mathbb{S}^{n})$, $D<\pi$ then ${\rm diam}(M)\geq D$. Thus if $K_{M}\geq 1$ and $\lambda_{1,p}(M) \leq C(n,p,\pi/2) \cdot \lambda_{1,p}(\mathbb{S}^{n})$ then by the previous argument, we have ${\rm diam}(M)\geq \pi/2$. Hence, by the Grove-Shiohama diameter sphere theorem, $M$ is homeomorphic to $\mathbb{S}^n$, completing the proof of Theorem \ref{thm2}. If ${\rm Ric}_{M}\geq (n-1)$, ${\rm inj}_{M}\geq i_0>0$ and $\lambda_{1,p}(M) \leq C(n,p, \pi-\epsilon(n,i_0)) \cdot \lambda_{1,p}(\mathbb{S}^{n})$ then  ${\rm diam}(M)\geq \pi-\epsilon(n,i_0)$. As a consequence of Bessa's diameter sphere theorem, the manifold \( M \) is diffeomorphic to the standard sphere \( \mathbb{S}^n \), which completes the proof of Theorem \ref{thm1}.

\section{Acknowledgments}
    I would like to express my gratitude to G. Pacelli Bessa for insightful discussions that significantly enriched this work. I also acknowledge that Paulo Silva was partially supported by CAPES‑Brazil under Grant No. 88881.126989/2025‑01.

\end{document}